\documentclass[12pt]{article}
\usepackage[a4paper, total={6in, 8in}]{geometry}
\usepackage[utf8]{inputenc} 
\usepackage[T1]{fontenc} 
\usepackage{amsmath}
\usepackage{amsthm}
\usepackage{setspace}
\usepackage{mathpazo} 
\usepackage{amssymb}
\usepackage{float}
\usepackage{graphicx}
\usepackage{xcolor}
\usepackage[maxbibnames=99]{biblatex}
\usepackage{ytableau}
\usepackage{hyperref}
\usepackage[autostyle=true]{csquotes}

\newtheorem{teo}{Theorem}
\newtheorem*{teo*}{Theorem}

\newtheorem{lemma}[teo]{Lemma}

\newtheorem*{conjecture*}{Conjecture}

\newtheorem{prop}[teo]{Proposition}
\theoremstyle{remark}
\newtheorem{rmk}[teo]{Remark}
\theoremstyle{definition}

\theoremstyle{definition}
\newtheorem{mydef}[teo]{Definition}
\newtheorem{esempio}[teo]{Example}

\begin{document}
\begin{center}
\textbf{\Large A deletion-contraction long exact sequence for chromatic symmetric homology} 
\vspace{5mm}
 \\Azzurra Ciliberti\footnote{La Sapienza Università di Roma - azzurra.ciliberti@uniroma1.it} 
 
 \end{center}
 \begin{abstract}
  \noindent 
  In \cite{MR4093019}, the authors generalize Stanley's chromatic symmetric function \cite{3} to vertex-weighted graphs. In this paper we find a categorification of their new invariant extending the definition of chromatic symmetric homology to vertex-weighted graphs. We prove the existence of a deletion-contraction long exact sequence for chromatic symmetric homology which gives a useful computational tool and allow us to answer two questions left open in \cite{2}. {In particular, we prove that, for a graph $G$ with $n$ vertices, the maximal index with nonzero homology is not greater that $n-1$. Moreover, we show that the homology is non-trivial for all the indices between the minimum and the maximum with this property.}  
\end{abstract}
\section*{Introduction}
The \emph{chromatic symmetric function} $X_G$ of a graph $G$, defined by Stanley in \cite{3}, is a remarkable combinatorial invariant which refines the chromatic polynomial. In \cite{1}, Sazdanovic and Yip categorified this invariant by defining a new homological theory, called the \emph{chromatic symmetric homology} of $G$. This construction, inspired by Khovanov's categorification of the Jones polynomial \cite{5}, is obtained by assigning a graded representation of the symmetric group to every subgraph of $G$, and a differential to every cover relation in the Boolean poset of subgraphs of $G$. The chromatic symmetric homology $H_{\ast,\ast}(G)$ is then defined as the homology of this chain complex; its bigraded Frobenius series $Frob_G(q,t)$, when evaluated at $q = t = 1$, reduces
to Stanley's chromatic symmetric function expressed in the Schur basis.
This categorification has interesting properties which {have been} investigated in \cite{2} and \cite{MR4540909}.

In \cite{MR4093019}, Logan Crew and Sophie Spirkl generalize Stanley's chromatic symmetric function \cite{3} to vertex-weighted graphs $(G,w)$ with the definition of the $weighted$ $chromatic$ $symmetric$ $function$ $X_{(G,w)}$. One of the primary motivations for extending the chromatic symmetric function to vertex-weighted graphs is the existence of a deletion-contraction relation in this setting, which, as known, holds for the chromatic polynomial, but doesn't hold for the chromatic symmetric function, as observed by Stanley in \cite{3}.

In this paper we generalize chromatic symmetric homology to vertex-weighted graphs. We obtain in this way a categorification of the weighted chromatic symmetric function that we call $weighted$ $chromatic$ $symmetric$ $homology$ and we denote by $H_{\ast,\ast}(G,w)$.
{The weighted chromatic symmetric homology specializes to the the chromatic symmetric homology if $w=\bold{1}$ is the function assigning weight 1 to each vertex, i.e. if $G$ is an unweighted graph.}

Moreover, we prove the existence of a deletion-contraction long exact sequence for the weighted chromatic symmetric homology which lifts to homology the 

deletion-contraction relation that holds for the function defined by Crew and Spirkl. 

{In particular}, we prove that

\begin{teo*}
Let $(G,w)$ be a vertex-weighted graph and let $e$ be an edge of $G$.
For each $j \geq 0$, there is a long exact sequence in homology 
\begin{center}
    $\to H_{i,j}(G\setminus e,w) \to H_{i,j}(G,w) \to  H_{i-1,j}(G/e,w/e) \xrightarrow{} H_{i-1,j}(G\setminus e,w) \to \dots$,
    \end{center}
where $G\setminus e$ denotes the graph $G$ with the edge $e$ removed, $G/e$ denotes the graph $G$ with the edge $e$ contracted to a point, {and $w/e$ denotes the weight function on $G/e$ defined in Section \ref{s1}}.
\end{teo*}

The long exact sequence in homology gives a useful computational tool and allow us to answer two questions left open in \cite{2}.

Let $span_0(G)$ denote the homological span of the degree 0 chromatic symmetric homology of $G$.  In \cite{2}, the authors formulate the following two conjectures.

\begin{conjecture*}[C.5]
Given any graph $G$, chromatic symmetric homology groups $H_{i,0}(G;\mathbb{C})$ are non-trivial for all $0 \leq i \leq span_0(G) - 1$.
\end{conjecture*}
\begin{conjecture*}[C.6]
Let $G$ be a graph with $n$ vertices and $m$ edges, and let $b$ denote the number of blocks of $G$. Then $n-b\leq span_0(G)\leq n-1$.
\end{conjecture*}

Using the deletion-contraction long exact sequence for chromatic symmetric homology {we show that Conjecture C.5 and a part of Conjecture C.6 are true}, also for the case of vertex-weighted graphs.

{
In particular, denoting by $k_{max}^j(G,w)$ the largest index $k$ such that $H_{k,j}(G,w)\neq 0$ and by $k_{min}^j(G,w)$ the smallest one ($k_{min}^0(G,w)$ is always 0 in the case of simple graphs), we prove that}

{
\begin{teo*}
Given any graph $(G,w)$, chromatic symmetric homology groups $H_{i,j}(G,w;\mathbb{C})$ are non-trivial for all $k_{min}^j(G,w) \leq i \leq k_{max}^j(G,w)$, $j\geq 0$. 
\end{teo*}}
{
\begin{teo*}
Let $(G,w)$ be a graph with $n$ vertices and $m$ edges. Then $k_{max}^j(G,w)\leq n-1$ for all $j\geq 0$. Moreover, if $m\geq 1$, $k_{max}^0(G,w)\leq n-2$, so $span_0(G)\leq n-1$.
\end{teo*}}

The paper is organized as follows. In Section 1 we recall the definition and some basic properties of the weighted chromatic symmetric function. In Section 2 we build our categorification and prove the existence of a long exact sequence in homology that lifts the deletion-contraction relation for the weighted chromatic symmetric function. Finally, in Section 3, we present some applications of the mentioned sequence and we prove the {last two theorems} above.

\section{Weighted chromatic symmetric function}\label{s1}
{Let $G$ be a graph. Then $G\setminus e$ denotes the graph $G$ with the edge $e$ removed and $G/e$ denotes the graph $G$ with the edge $e$ contracted to a point.} 
\begin{mydef}
Define a $vertex$-$weighted$ $graph$ $(G,w)$ to be a graph $G=(V(G), E(G))$ together with a vertex-weight function
$w : V (G) \to \mathbb{N}$. The $weight$ of a vertex $v \in V (G)$ is $w(v)$.
\end{mydef}
\begin{rmk}
Let $G$ be any graph. Then $G$ can be viewed as the vertex-weighted graph $(G,\bold{1})$, where $\bold{1}$ is the function assigning weight 1 to each vertex.
\end{rmk}
\begin{mydef}
 Given a vertex-weighted graph $(G, w)$, {we say that $F \subseteq V (G)$ is a $state$ $of$ $G$}, and we define the $total$ $weight$ $w(F)$ of $F$ to be $\displaystyle\sum_{v \in F} w(v)$. Moreover, we define the total weight $w(G)$ of $G$ to be the total weight of $V(G)$.
\end{mydef}

{The set $Q(G)$ of all the states of $G$ has a stucture of Boolean lattice, ordered by reverse inclusion. In the Hasse diagram of $Q(G)$, we direct an edge $\epsilon (F,F')$ from a subgraph $F$ to a  subgraph $F'$ if and only if $F'$ can be obtained by removing an edge from $F$.} 

In \cite{MR4093019}, Logan Crew and Sophie Spirkl generalize Stanley's chromatic symmetric function \cite{3} to vertex-weighted graphs with the following definition:

\begin{mydef}
 Let $(G,w)$ be a vertex-weighted graph. Then the $weighted$ $chromatic$ $symmetric$ $function$ is 
 \begin{center}
  $X_{(G,w)}(x_1,x_2,\dots) = \displaystyle\sum_\kappa \prod_{v \in V(G)}x_{\kappa(v)}^{w(v)}$,   
 \end{center}
 where the sum ranges over all proper colorings $\kappa : V(G) \to \mathbb{N}$ of $G$.
\end{mydef}
\begin{rmk}
If $G$ has a loop, then $X_{(G,w)}=0$ for every $w : V (G) \to \mathbb{N}$. Moreover, if $e_1$, $e_2$ are edges of $G$ with the same endpoints, then  $X_{(G,w)}=X_{(G\setminus e_1,w)}=X_{(G\setminus e_2,w)}$ for every $w : V (G) \to \mathbb{N}$.
\end{rmk}
\begin{rmk}
 {Note that $X_{(G,\bold{1})}=X_G$, where $X_G$ is the usual chromatic symmetric function.}
\end{rmk}

Recall that, if $\lambda=(\lambda_1, \dots, \lambda_k)$ is a partition of a positive integer $n$, i.e. a non-increasing sequence of positive integers whose sum is $n$, the power sum symmetric function $p_\lambda$ is defined as
\begin{center}
    $p_\lambda (x_1, x_2, \dots) = p_{\lambda_1}(x_1, x_2, \dots)\cdots p_{\lambda_k}(x_1, x_2, \dots)$, 
\end{center}
where $p_r(x_1, x_2, \dots) = x_1^r + x_2^r + \dots$, for $r \in \mathbb{N}$.

{Let $\Lambda_n$ be the $\mathbb{Z}$-module of the homogeneous symmetric functions of degree $n$. Then $\{$ $p_\lambda$ $|$ $\lambda$ partition of $n$ $\}$ is a basis for $\Lambda_n$. Another basis for $\Lambda_n$ is given by the Schur symmetric functions $\{$ $s_\lambda$ $|$ $\lambda$ partition of $n$ $\}$.
Moreover, let {$\Lambda^\mathbb{C}=\displaystyle\bigoplus_{n\geq 0}\Lambda_n$} denote the space of symmetric functions in the indeterminates $x_1, x_2, \dots$.}

\begin{mydef}
 Given a vertex-weighted graph $(G,w)$, and $F \subseteq E(G)$, we define $\lambda(G,w,F)$ to be the partition of $w(G)$
whose parts are the total weights of the connected components of $(G', w)$, where $G' = (V(G),F)$.
\end{mydef}

\begin{lemma}[{\cite{MR4093019}, Lemma 3}]\label{l1}
Let $(G,w)$ be a vertex-weighted graph. Then
\begin{center}
$X_{(G,w)} = \displaystyle\sum_{F \subseteq E(G)} (-1)^{|F|} p_{\lambda(G,w,F)}$.
\end{center}
\end{lemma}

One of the primary motivations for extending the chromatic symmetric function to vertex-weighted graphs is the existence of a deletion-contraction relation in this setting. 
{
\begin{mydef}
  Let $(G,w)$ be a vertex-weighted graph, and let $e = (v_1,v_2) \in E(G)$. We define $w/e : V(G/e) \to \mathbb{N}$ to be the modified weight function on $G/e$ such that $w/e = w$ if $e$ is a loop, and otherwise $(w/e)(v) = w(v)$ if $v \neq v_1,v_2$, and for the vertex $v^\ast$ of $G/e$ formed by the contraction, $(w/e)(v^\ast) = w(v_1) + w(v_2)$.   
\end{mydef}
}
We have the following:

\begin{teo}[{\cite{MR4093019}, Lemma 2}]\label{t1}
 Let $(G, w)$ be a vertex-weighted graph, and let $e \in E(G)$ be any edge. Then
 \begin{center}

$X_{(G,w)} = X_{(G \setminus e,w)} - X_{(G/e,w/e)}$.
     
 \end{center}
\end{teo}

{Note that the deletion-contraction relation of Theorem \ref{t1} does not give a similar relation for the ordinary chromatic symmetric function, since if we contract a non loop edge we do not get an ordinary chromatic symmetric function.}

\section{Weighted chromatic symmetric homology}
Now we build a categorification of the invariant just introduced.

{In this section we assume that the set of edges of $G$ is ordered.} 

Let $\frak S_n$ denote the symmetric group on $n$ elements. The irreducible representations of $\frak S_n$ over $\mathbb{C}$ are indexed by the partitions of $n$, and are called $Specht$ $modules$. Let $\bold{S}^\lambda$ denote the Specht module indexed by $\lambda$. 

The Grothendieck group $R_n$ of representations of $\frak S_n$ is the free abelian group on the isomorphism classes $[\bold{S}^\lambda]$ of irreducible representations of $\frak S_n$, modulo the subgroup generated by all $[V \oplus W]-[V]-[W]$. Let $R = \displaystyle\bigoplus_{n \geq 0 } R_n$. If $[V]\in R_a$ and $[W]\in R_b$, {define a multiplication in $R$ by}
 \begin{center}
     $[V]\circ [W]=[Ind_{\frak S_a \times \frak S_b}^{\frak S_{a+b}}V \otimes W]$.
 \end{center}
{Here the tensor product $V \otimes W$ is regarded as a representation of $\frak S_n \times \frak S_m$ in the obvious way: $(\sigma \times \tau) \cdot (v \otimes w) = \sigma \cdot v \otimes \tau \cdot w$; and $\frak S_n \times \frak S_m$ is regarded as a subgroup of $\frak S_{n+m}$ with $\frak S_n$ acting on the first $n$ integers and $\frak S_m$ acting on the last $m$ integers. The induced representation can be defined quickly by the formula
\begin{center}
    $Ind_{\frak S_n \times \frak S_m}^{\frak S_{n+m}}=\mathbb{C}[\frak S_{n+m}]\otimes_{\mathbb{C}[\frak S_n \times \frak S_m]}(V \otimes W)$.
\end{center}
It is straightforward to verify that this product is well defined and makes R into a commutative, associative, graded ring with unit.}

 The morphism of graded rings given by sending the Specht modules to the Schur functions

 \begin{center}
   $ch : R \to \Lambda^\mathbb{C}$, $[\bold{S}^\lambda]\to s_\lambda$  
 \end{center}

 is an isomorphism. 
 
Moreover, for $n \in \mathbb{N}$, we have

\begin{equation}\label{eqn}
ch^{-1}(p_n)=\displaystyle\sum_{i=0}^{n-1}(-1)^i [\bold{S}^{(n-i,1^i)}].
\end{equation}
    
{For the proofs of these two last facts see \cite{4}, Section 7.3.}

With the notation of \cite{1}, we define:

\begin{mydef}
Let $(G,w)$ be a vertex-weighted graph. Suppose $F \subseteq E(G)$ is a state with $r$ connected components of total weights $b_1^w,\dots,b_r^w$ respectively. To $F$, we assign the {graded} $\frak S_{w(G)}$-module
\begin{equation}
    M_F^w = Ind_{\frak S_{b_1^w} \times \dots \times \frak S_{b_r^w}} ^{\frak S_{w(G)}}(\bold{L}_{b_1^w}\otimes \dots \otimes \bold{L}_{b_r^w} ),
\end{equation}
where $\bold{L}_a$ denotes the $q$-graded $\frak S_a$-module
\begin{equation}
    \bold{L}_a = \displaystyle\bigoplus_{j=0}^{a-1}\bold{S}^{(a-j,1^j)},
\end{equation}
and $\bold{S}^{(a-j,1^j)}$ is the Specht module related to the partition $(a-j,1^j)$ of the positive integer $a$. The {grading is given by the index $j$}.
\end{mydef}

\begin{mydef}
For $i\geq 0$, the $i$-th $weighted$ $chain$ $module$ for $(G,w)$ is
\begin{center}
    $C_i(G,w)=\displaystyle\bigoplus_{|F|=i}M_F^w$.
\end{center}
More precisely, since $M_F^w = \displaystyle\bigoplus_{j\geq 0}(M_F^w)_j$ is graded, then for $i,j \geq 0$, we define
\begin{center}
    $C_{i,j}(G,w)=\displaystyle\bigoplus_{|F|=i}(M_F^w)_j$.
\end{center}
\end{mydef}

\begin{rmk}
Observe that $(M_F^w)_j = 0$ if $j \geq b_t^w$ for all $t=1,\dots,r$.
\end{rmk}

Since the differential defined in \cite{1} depends only on the $b_i$'s, we can define a differential in the same way, replacing the $b_i$'s with the $b_i^w$'s.

Let $F$ be a state of $G$. Suppose $F' = F - e$ where $e \in E(G)$. We define the $\frak S_{w(G)}$-modules morphism $d_\epsilon^{(G,w)} : M_F^w \to M_{F'}^w$, i.e. the $per$-$edge$ $maps$, in the following way.

 There are two cases to consider:
\begin{itemize}
    \item [Case 1] The edge $e$ is incident to vertices in the same connected component of $F'$. Since $M_F^w$ and $M_{F'}^w$ are {equal}, we define $d_\epsilon : M_F^w \to M_{F'}^w$ to be the identity map.
    \item [Case 2] The edge $e$ is incident to vertices in different connected components of $F'$. First, consider the simplest case where $F$ consists of one connected component and $F'$ consists of two components $A$ and $B$. Suppose $w(A) = a$ and $w(B)=b$, so that $a+b=w(G)$. 
    Since, by Frobenius Reciprocity, $Hom_{\frak S_{w(G)}}(M_F^w,M_{F'}^w) \cong Hom_{\frak S_a \times \frak S_b} (\Lambda^\ast T \oplus (\Lambda^\ast T)[1], \Lambda^\ast T )$, where  
    
    $T = ( \bold{S}^{(a-1,1)} \otimes \mathbb{1}_{\frak S_b}) \oplus (\mathbb{1}_{\frak S_a} \otimes \bold{S}^{(b-1,1)})$ (see \cite{1}, Lemma 2.6), we choose the element $d_\epsilon \in Hom_{\frak S_{w(G)}}(M_F^w,M_{F'}^w)$ to be the map that corresponds to the $(\frak S_a \times \frak S_b)$-module map that is the identity on $ \Lambda^\ast T$ and zero on $(\Lambda^\ast T)[1]$. In the general case when $F$ has more than one connected component, the definition of the per-edge map is achieved {by recursion} on the two-component case.

    {
    Suppose $F$ is a state with $r$ connected components $B_1, \dots, B_r$ of total weights $b_1^w,\dots,b_r^w$. Further suppose that the removal of the edge $e \in E(G)$ decomposes $B_r$ into two components $A$ and $B$ of total weights $a$ and $b$ respectively ($a+b=b_r^w$).}
 {  
Let $d_\zeta : \bold{L_{b_r^w}} \to Ind_{\frak S_a \times \frak S_b}^{\frak S_{b_r^w}} (\bold{L_a}\otimes \bold{L_b})$ be the per-edge map defined previously (note that $M_{B_r}^{b_r^w}=\bold{L}_{b_r^w}$, since $B_r$ is connected), and let $\bold{N}=\bold{L}_{b_1^w}\otimes \cdots \otimes \bold{L}_{b_{r-1}^w}$.}
{
The map $d_\epsilon : M_F^w  \to M_{F'}^w$ is chosen to be 
\begin{center}
    $d_\epsilon = Ind_{\frak S_{b_1^w} \times \dots \times \frak S_{b_{r-1}^w} \times \frak S_{b_r^w}}^{\frak S_{w(G)}}(id_{\bold{N}}\otimes d_\zeta)$
\end{center}}

    \end{itemize}

{
\begin{mydef}
    Let $F$ and $F'$ be states of $G$. Assume that $F'=F\setminus e$, $e \in E(F)$. The sign of $\epsilon=\epsilon(F,F')$, $sgn(\epsilon)$, is defined as $(-1)^k$, where $k$ is the number of edges of $F$ less than $e$.
\end{mydef}
}
    
\begin{mydef}
 For $i\geq 0$, define $d_i^{(G,w)} : C_i(G,w) \to C_{i-1}(G,w)$ letting
 \begin{center}
     $d_i^{(G,w)} = \displaystyle\sum_\epsilon sgn(\epsilon)d_\epsilon^{(G,w)}$,
 \end{center}
 where the sum is over all edges $\epsilon$ in the Hasse diagram of $Q(G)$ joining a state with $i$ edges to a state with $i-1$ edges. We also define $d_{i,j}^{(G,w)} : C_{i,j}(G,w) \to C_{i-1,j}(G,w)$ to be the map $d_i^{(G,w)}$ in the $j$-th grading. 
 \end{mydef}

 \begin{prop}
 {The maps $d_i^{(G,w)}$ form a differential on the chain complex $C_\ast(G, w)$.}
 \end{prop}
 
 \begin{proof}
 The proof is completely analogous to that of Proposition 2.10 of \cite{1} replacing the $b_i$'s with the $b_i^w$'s.
 \end{proof}
 
 \begin{mydef}
For $i,j \geq 0$, the $(i,j)$-th $weighted$ $chromatic$ $symmetric$ $homology$ of $(G,w)$ is
\begin{center}
    $H_{i,j}(G,w) = \mathrm{ker}\, d_{i,j}^{(G,w)}/\mathrm{im}\, d_{i+1,j}^{(G,w)}$.
\end{center}
{Moreover, we define}
\begin{center}
    $H_i(G,w) = \displaystyle\bigoplus_{j \geq 0} H_{i,j}(G,w)$.
\end{center}
\end{mydef}

\begin{rmk}
    {$H_{\ast,\ast}(G,\bold{1})=H_{\ast,\ast}(G)$, where $H_{\ast,\ast}(G)$ is the usual chromatic symmetric homology.}
\end{rmk}

\begin{esempio}\label{ex1}
Let $(K_2,w)$ be the segment with a vertex $v_1$ of weight 1 and the other $v_2$ of weight 2. The labels of the vertices indicate their weights.
\begin{figure}[H]
\centering
\includegraphics[width= 0.7 \textwidth]{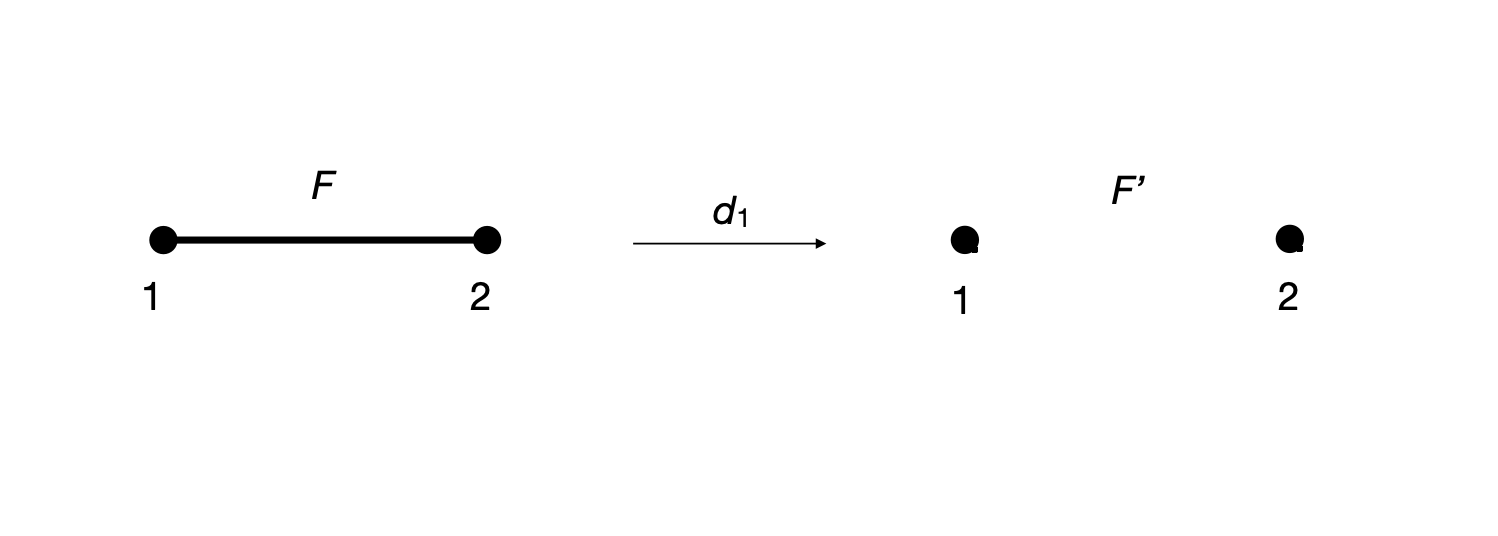}
\label{f4}
\end{figure}
We have 
\begin{itemize}
    \item $C_{1,0}(K_2,w)= (M_F^w)_0=\bold{S}^{(3)}$;
    \item $C_{0,0}(K_2,w)= (M_{F'}^w)_0=Ind_{\frak S_2\times \frak S_1}^{\frak S_3}\bold{S}^{(2)}\otimes \bold{S}^{(1)}=\bold{S}^{(3)}\oplus \bold{S}^{(2,1)}$;
    \item $C_{1,1}(K_2,w)=(M_F^w)_1=\bold{S}^{(2,1)}$;
    \item $C_{0,1}(K_2,w)= (M_{F'}^w)_1=Ind_{\frak S_2\times \frak S_1}^{\frak S_3}\bold{S}^{(1,1)}\otimes \bold{S}^{(1)}=\bold{S}^{(2,1)}\oplus \bold{S}^{(1^3)}$;
    \item $C_{1,2}(K_2,w)= (M_F^w)_2=\bold{S}^{(1^3)}$;
    \item $C_{0,2}(K_2,w)= 0$.
\end{itemize}
Therefore, $H_{1,0}(K_2,w)=H_{1,1}(K_2,w)=H_{0,2}(K_2,w)=0$, $H_{0,0}(K_2,w)=\bold{S}^{(2,1)}$, $H_{0,1}(K_2,w)=H_{1,2}(K_2,w)=\bold{S}^{(1^3)}$.

In general, 
\begin{itemize}
    \item $C_{1,0}(K_2,w)= (M_F^w)_0=\bold{S}^{(w(v_1)+w(v_2))}$;
    \item $C_{0,0}(K_2,w)= (M_{F'}^w)_0=Ind_{\frak S_{w(v_1)}\times \frak S_{w(v_2)}}^{\frak S_{w(v_1)+w(v_2)}}\bold{S}^{(w(v_2))}\otimes \bold{S}^{(w(v_2))}$
    
    \hspace{4cm}$=\bold{S}^{(w(v_1)+w(v_2))}\oplus \displaystyle\bigoplus_{\lambda}(\bold{S}^\lambda)^{m_\lambda}$.
\end{itemize}

{
We don't give the details about the $\bold{S}^\lambda$'s which appear in the last formula and their multiplicities. You can find an explanation of it in \cite{4}, Section 7.3. We say only that they are all different from $\bold{S}^{(w(v_1)+w(v_2))}$. Therefore, we have $H_{1,0}(K_2,w)=0$ and $H_{0,0}(K_2,w) \neq 0$. Moreover, $H_{i,0}(K_2,w)=0$ for any $i\geq 2$, since $K_2$ does not have any states with more than one edge.}

\end{esempio}

\begin{mydef}

The bigraded $Frobenius$ $series$ of {$H_{\ast,\ast} (G,w) = \displaystyle\bigoplus_{i,j \geq 0} H_{i,j}(G,w)$} is
\begin{center}
    $Frob_{(G,w)}(q,t) = \displaystyle\sum_{i,j \geq 0} (-1)^{i+j} t^i q^j ch(H_{i,j}(G,w))$.
\end{center}
\end{mydef}

\begin{esempio}
Let's consider the vertex-weighted graph of the previous example. We have
\begin{center}
    $Frob_{(K_2,w)}(q,t)= - (q+tq^2)s_{(1^3)}+ s_{(2,1)}$.
\end{center}
\end{esempio}

\begin{lemma}\label{l2}
For any vertex-weighted graph $(G,w)$,
\begin{center}
    $\displaystyle\sum_{i,j \geq 0} (-1)^{i+j} ch(H_{i,j}(G,w)) = \displaystyle\sum_{i,j \geq 0} (-1)^{i+j} ch(C_{i,j}(G,w))$.
\end{center}
\end{lemma}
\begin{proof}
 Let $n$ be any positive integer. Any short exact sequence of $\frak S_n$-modules $0 \to A \to B \to C \to 0$ is split exact, so $B \cong A \oplus C$ and $ch(B) = ch(A)+ch(C)$.\\
Let $Z_{i,j}(G,w) = \mathrm{ker}\,d_{i,j}^{(G,w)}$ and $B_{i,j}(G,w) = \mathrm{im}\,d_{i+1,j}^{(G,w)}$. For $i,j \geq 0$, we have short exact sequence $0 \to Z_{i,j}(G,w) \to C_{i,j}(G,w) \to B_{i-1,j}(G,w) \to 0$ and 
$0 \to B_{i,j}(G,w) \to Z_{i,j}(G,w) \to H_{i,j}(G,w) \to 0$, where $B_{-1,j}(G,w)$ is understood to be zero. Thus

\vspace{0.5cm}
$ch( C_{i,j}(G,w)) = ch(Z_{i,j}(G,w)) + ch (B_{i-1,j}(G,w)) $
\vspace{0.3cm}
    
    \hspace{2.6cm}$= ch( H_{i,j}(G,w)) + ch( B_{i,j}(G,w)) + ch( B_{i-1,j}(G,w)) $.

\vspace{0.5cm}
{If we multiply this by $(-1)^{i+j}$ and we sum over all $i,j\geq 0$, we get:}
\vspace{0.5cm}

\hspace{-0.7cm}{$\displaystyle\sum_{i,j \geq 0} (-1)^{i+j} ch(C_{i,j}(G,w)) =\displaystyle\sum_{i,j \geq 0} (-1)^{i+j} ch(H_{i,j}(G,w))+\displaystyle\sum_{i,j \geq 0} (-1)^{i+j} ch(B_{i,j}(G,w))+$
\hspace{-0.7cm}$\displaystyle\sum_{i,j \geq 0} (-1)^{i+j} ch(B_{i-1,j}(G,w)) = \displaystyle\sum_{i,j \geq 0} (-1)^{i+j} ch(H_{i,j}(G,w))+\displaystyle\sum_{i,j \geq 0} (-1)^{i+j} ch(B_{i,j}(G,w))$} 

\hspace{-0.6cm}{$-\displaystyle\sum_{t,j \geq 0} (-1)^{t+j} ch(B_{t,j}(G,w))= \displaystyle\sum_{i,j \geq 0} (-1)^{i+j} ch(H_{i,j}(G,w))$. }
\end{proof}



\begin{teo}\label{t3}
Weighted chromatic symmetric homology categorifies the weighted chromatic symmetric
function. That is, for any vertex-weighted graph $(G,w)$,
\begin{center}
    $Frob_{(G,w)}(1,1) = X_{(G,w)}$.
\end{center}
\end{teo}
\begin{proof}
Using Lemma \ref{l2}, \ref{eqn} and Lemma \ref{l1}, we have\\

$\hspace{-0.5cm}Frob_{(G,w)}(1,1)   = \displaystyle\sum_{i,j \geq 0} (-1)^{i+j} ch(H_{i,j}(G,w)) 
= \displaystyle\sum_{i\geq 0}(-1)^i \Big(\displaystyle\sum_{j\geq 0}(-1)^j ch(C_{i,j}(G,w))\Big)$ 

\hspace{1.5cm}$= \displaystyle\sum_{i\geq 0} (-1)^i \displaystyle\displaystyle\sum_{F \subseteq E(G) : |F| = i} p_{\lambda(G,w,F)} = X_{(G,w)}$.
\end{proof}

Now we want to lift to homology the result of Theorem \ref{t1}.

\begin{prop}
\label{t2}
Let $(G,w)$ be a vertex-weighted graph and let $e$ be an edge of $G$.For each $i,j \geq 0$, there is a short exact sequence of $\frak S_{w(G)}$-modules
\begin{center}
    $0\to C_{i,j}(G\setminus e,w) \to C_{i,j}(G,w) \to C_{i-1,j}(G/e,w/e) \to 0$.
\end{center}
\end{prop}

\begin{proof}
{By definition} 

\begin{center}
{    
$\displaystyle C_{i,j}(G\displaystyle\setminus e,w)=\displaystyle\bigoplus_{\displaystyle|F\displaystyle|=i, F \displaystyle\subseteq E(G \displaystyle\setminus e)}(M_F^w)_j$}

{
\hspace{3cm}$=\displaystyle\bigoplus_{\displaystyle|F\displaystyle|=i, F \displaystyle\subseteq E(G), e \displaystyle\notin F}\displaystyle(M_F^w)_j$} 

{
\hspace{4.3cm}$\displaystyle\subseteq \displaystyle\bigoplus_{\displaystyle|F\displaystyle|=i, F \displaystyle\subseteq E(G)}(M_F^w)_j$
$= C_{i,j}(G,w)$. }
\end{center}
Therefore, there is a short exact sequence
\begin{center}
    $0\to C_{i,j}(G\setminus e,w) \xrightarrow{\iota_i} C_{i,j}(G,w) \xrightarrow{\pi_i} \frac{C_{i,j}(G,w)}{C_{i,j}(G\setminus e,w)} \to 0$,
\end{center}
where $\iota_i$ is the inclusion and $\pi_i$ is the projection to the quotient. 

We have that
\begin{center}
  $\frac{C_{i,j}(G,w)}{C_{i,j}(G\setminus e,w)}=\frac{\displaystyle\bigoplus_{|F|=i, F \subseteq E(G)}(M_F^w)_j}{\displaystyle\bigoplus_{|F|=i, F \subseteq E(G), e \notin F}(M_F^w)_j} \cong \displaystyle\bigoplus_{|F|=i, F \subseteq E(G), e \in F}(M_F^w)_j$.  
\end{center}

Since, if $F$ is a state of $(G,w)$ with $i$ edges such that $e \in F$, then $M_F^w=M_{F/e}^{w/e}$, because the contraction does not change the total weight of the connected components of $F$, and $F/e$ is a state of $(G/e,w/e)$ with $i-1$ edges, we have that
\begin{center}
  $\displaystyle\bigoplus_{|F|=i, F \subseteq E(G), e \in F}(M_F^w)_j=C_{i-1,j}(G/e,w/e)$,  
\end{center}
and the theorem follows.
\end{proof}

\begin{rmk}\label{t_2}
If $G$ is an unweighted graph, for each $i,j \geq 0$, we have the following short exact sequence of $\frak S_{|V(G)|}$-modules
\begin{center}
    $0\to C_{i,j}(G\setminus e) \to C_{i,j}(G) \to C_{i-1,j}(G/e,\bold{1}/e) \to 0$.
\end{center}
\end{rmk}

\begin{prop}\label{p1}
Let $(G,w)$ be a vertex-weighted graph and let $e$ be an edge of $G$.
For each $j \geq 0$, there is a short exact sequence of chain complexes
\begin{center}
    {$0\to C_{\ast,j}(G\setminus e,w) \to C_{\ast,j}(G,w) \to C_{\ast -1,j}(G/e,w/e) \to 0$.}
\end{center}
\end{prop}

\begin{proof}
With the notation of the proof of Proposition \ref{t2}, we have to show that, for each $i \geq 0$, $d_i^{(G,w)} \circ \iota_i = \iota_{i-1} \circ d_i^{(G \setminus e,w)}$ and $d_{i-1}^{(G/e,w/e)} \circ \pi_i = \pi_{i-1} \circ d_i^{(G,w)}$. It is clear that the first equality holds. 
Let's look at the second.

If $i=0,1$, we have $0$ on both sides. Consider $i\geq 2$. Since, if $F$ is a state of $(G,w)$ with $i$ edges such that $e \in F$, then $M_F^w=M_{F/e}^{w/e}$, $\pi_i$ is the map such that
\begin{center}
    
   { $\pi_{i_{|_{ M_F^w}}} = \begin{cases}
\text{id \hspace{0.8cm} if $e \in F$,}\\
 \text{0  \hspace{1cm} if $e \notin F$.}
\end{cases}$}
    
\end{center}
Therefore, 

$\pi_{i-1} \circ d_i^{(G,w)}= \displaystyle\sum_{\epsilon} sgn(\epsilon)\pi_{i-1}\circ d_\epsilon^{(G,w)}=\displaystyle\sum_{\epsilon'} sgn(\epsilon')d_{\epsilon'}^{(G,w)}$, where the last sum is over all the $\epsilon'$ in the Hasse diagram of $Q(G,w)$ joining a state of $(G,w)$ with $i$ edges that contains $e$ to a state of $(G,w)$ with $i-1$ edges that also contains $e$.

On the other hand, $d_{i-1}^{(G/e,w/e)} \circ \pi_i= \displaystyle\sum_{\epsilon''} sgn(\epsilon'')d_{\epsilon''}^{(G/e,w/e)}$, where the sum is over all the $\epsilon''$ in the Hasse diagram of $Q(G/e,w/e)$ joining a state of $(G/e,w/e)$ with $i-1$ edges to a state of $(G/e,w/e)$ with $i-2$ edges.

We know that, if $F$ is a state of $G$ with $i$ edges such that $e \in F$, then $M_F^w=M_{F/e}^{w/e}$ and $F/e$ is a state of $(G/e,w/e)$ with $i-1$ edges. Therefore, if $\epsilon'$ is an edge in the Hasse diagram of $Q(G,w)$ connecting a state $F$ of $(G,w)$ with $i$ edges that contains $e$ with a state $F'$ of $(G,w)$ with $i-1$ edges that also {contains} $e$,
\begin{center}
    $d_{\epsilon'}^{(G,w)} : M_F^w=M_{F/e}^{w/e} \to M_{F'}^w=M_{F'/e}^{w/e}$
\end{center}
coincides with $d_{\epsilon''}^{(G/e, w/e)}$, where $\epsilon''$ is an edge in the Hasse diagram  of $Q(G/e,w/e)$ joining the state $F/e$ of $(G/e,w/e)$ with $i-1$ edges to the state $F'/e$ of $(G/e,w/e)$ with $i-2$ edges.

Since there is a bijection between the states of $G$ with $i$ edges that contains $e$ and the states of $(G/e,w/e)$ with $i-1$ edges, we have that the two sums coincide. Therefore,

\begin{center}
    $d_{i-1}^{(G/e,w/e)} \circ \pi_i = \pi_{i-1} \circ d_i^{(G,w)}$.
\end{center}
\end{proof}


Therefore, we have:

\begin{teo}\label{t4}
Let $(G,w)$ be a vertex-weighted graph and let $e$ be an edge of $G$.
For each $j \geq 0$, there is a long exact sequence in homology 
\begin{equation}\label{eqn}
\to H_{i,j}(G\setminus e,w) \to H_{i,j}(G,w) \to  H_{i-1,j}(G/e,w/e) \xrightarrow{\gamma^\ast} H_{i-1,j}(G\setminus e,w) \to \dots
\end{equation}
\end{teo}
\begin{proof}
The short exact sequences of chain complexes in Proposition \ref{p1} induce for each $j \geq 0$ a long exact sequence in homology.
\end{proof}
\begin{rmk}
The specialization of the Frobenius series at $q = t = 1$ recovers the deletion-contraction relation of Theorem \ref{t1}.
\end{rmk}
\begin{rmk}\label{r2}
The description for $\gamma^\ast$ follows from the standard diagram chasing argument in the zig-zag lemma and the result is as follows. It is the linear extension of the map that, given a state of $(G/e,w/e)$ with $i-1$ edges, where $e=(v_e,w_e)$ is an edge of $G$ {that has been contracted to a point}, expands $v_e=w_e$ by adding $e$ {with weight $w(v_e)$ at the vertex $v_e$ and $w(w_e)$ at the vertex $w_e$} and then deletes $e$. {In this way we get a state of $(G\setminus e, w)$ with $i-1$ edges.}   
\end{rmk}
\begin{rmk}
If $G$ is an unweighted graph, for each $j \geq 0$, we have the following long exact sequence in homology 
\begin{center}
    $\dots \longrightarrow H_{i,j}(G\setminus e) \rightarrow H_{i,j}(G) \rightarrow  H_{i-1,j}(G/e,\bold{1}/e) \xrightarrow{\gamma^\ast} H_{i-1,j}(G\setminus e) \longrightarrow \dots$.
\end{center}
\end{rmk}







\subsection{Properties of $H_{\ast,\ast}(G,w)$}

The deletion-contraction long exact sequence allows us to give a different and faster proof of the following two properties of chromatic symmetric homology, contained in \cite{1}, and to extend them to the case of vertex-weighted graphs.

\begin{prop}\label{p2}
If $(G,w)$ contains a loop, then $H_{\ast,\ast}(G,w)=0$.
\end{prop}
\begin{proof}
Let $(G,w)$ be a graph with a loop $l$. The exact sequence for $(G,w)$ with respect to $l$ is 
\begin{center}
  $\dots \to H_{i,j}(G/l,w/l) \xrightarrow{\gamma^\ast} H_{i,j}(G\setminus l,w) \rightarrow H_{i,j}(G,w) \rightarrow$  \\$H_{i-1,j}(G/l,w/l) \xrightarrow{\gamma^\ast} H_{i-1,j}(G\setminus l,w)\to \dots$.  
\end{center}
Using our description of the snake map $\gamma^\ast$ in Remark \ref{r2}, we get that the map
$H_{i,j}(G/l,w/l) \xrightarrow{\gamma^\ast} H_{i,j}(G\setminus l,w)$ is the identity map. Therefore, $H_{i,j} (G,w) = 0$ for all $i,j$.
\end{proof}

\begin{prop}
Let $(G,w)$ be a multigraph, i.e. a graph which is allowed to have multiple edges. Let $e_1$ and $e_2$ be two edges of $(G,w)$ with the same endpoints. Then $H_{\ast,\ast}(G,w) = H_{\ast,\ast}(G - e_2,w)$.
\end{prop}

\begin{proof}
In $G/e_2$ , $e_1$ becomes a loop so, by Proposition \ref{p2}, $H_{i,j}(G/e_2,w/e_2) = 0$ for all $i,j$. It follows from the long exact sequence {\ref{eqn}} that $H_{i,j}(G - e_2,w)$ and $H_{i,j}(G,w)$ are isomorphic modules. 
\end{proof}

{Therefore, from now on we assume that $G$ is simple, so without loops or multiple edges.}

Given two vertex-weighted graphs $(A,w_A)$ and $(B,w_B)$, let $(A+B,w_{A+B})$ denote their disjoint union, where 

$$
w_{A+B}(v)=\begin{cases}
			w_A(v), & \text{if $v \in V(A)$},\\
            w_B(v), & \text{if $v \in V(B)$}.
		 \end{cases}
$$

\begin{prop}
For $i,j \geq 0$,
\begin{center}
    $H_{i,j}(A+B,w_{A+B}) = \displaystyle\bigoplus_{\substack{p+r = i \\ q+s = j}} Ind_{\frak S_{w_A(A)} \times \frak S_{w_B(B)}}^{\frak S_{w_A(A)+w_B(B)}} (H_{p,q}(A,w_A)\otimes H_{r,s}(B,w_B))$.
\end{center}
\end{prop}
\begin{proof}
The proof is completely analogous to the unweighted case. See \cite{1}, Proposition 3.3.
\end{proof}

\begin{rmk}\label{r9}
If $(G,w)$ is a graph with homology $H_{i,j}(G,w)$ $ = \displaystyle\bigoplus_\lambda (\bold{S}^\lambda)^{\oplus m_\lambda}$, then the homology of the disjoint union of $G$ with a single vertex with weight $w_v$ is 
\begin{center}
    $H_{i,j}(G + \bullet) = \displaystyle\bigoplus_\mu (\bold{S}^\mu)^{\oplus m_\lambda}$,
\end{center}
where the sum is over all partitions $\mu$ which can be obtained by adding $w_v$ boxes to the partitions $\lambda$ indexing the irreducible factors of $H_{i,j}(G,w)$.
\end{rmk}

\section{Applications}

The deletion-contraction long exact sequence in homology has proved to be a useful computational tool. Moreover, we can use it to compute weighted chromatic symmetric homology starting from unweighted chromatic symmetric homology.
\begin{esempio}
Let $(K_2,w)$ be the segment with a vertex of weight 1 and the other of weight 2. We can compute its homology using the deletion-contraction long exact sequence. 

Let $G=P_3$ be the graph made of two segments with a vertex in common, and let $e \in E(G)$. We have that $(K_2,w)=G/e$ and $G \setminus e$ is the disjoint union of $K_2$ and an isolated vertex.

We have $H_{0,0}(G\setminus e)= H_{1,1}(G\setminus e)= \bold{S^{(2,1)}}\oplus \bold{S}^{(1^3)}$ and $H_{1,0}(G\setminus e)=0$.

Moreover, we have  $H_{0,0}(G)=H_{2,2}(G)=\bold{S}^{(1^3)}$, $H_{1,1}(G)=\bold{S}^{(2,1)}\oplus \bold{S}_{(1^3)}^{\oplus 2}$ and $H_{0,1}(G)=H_{2,0}(G)=H_{2,1}(G)=0$.

\vspace{0.5cm}
For $j=0$, we have the following long exact sequence in homology:

\begin{center}

$0 \longrightarrow H_{1,0}(K_2,w) \longrightarrow 0 \longrightarrow 0 \longrightarrow H_{0,0}(K_2,w) \longrightarrow $

$\longrightarrow \bold{S^{(2,1)}}\oplus \bold{S}^{(1^3)}\longrightarrow \bold{S}^{(1^3)} \longrightarrow 0$, 
\end{center}
from which we can conclude that $H_{1,0}(K_2,w) = 0$ and $H_{0,0}(K_2,w)= \bold{S}^{(2,1)}$.

\vspace{0.5cm}
For $j=1$, we have the following long exact sequence in homology:
\begin{center}

$0 \longrightarrow H_{1,1}(K_2,w) \longrightarrow \bold{S}^{(2,1)}\oplus \bold{S}^{(1^3)} \longrightarrow \bold{S}^{(2,1)}\oplus \bold{S}_{(1^3)}^{\oplus 2} \longrightarrow H_{0,1}(K_2,w) \longrightarrow 0$, 
\end{center}
from which we can conclude that $H_{1,1}(K_2,w) = 0$ and $H_{0,1}(K_2,w)= \bold{S}^{(1^3)}$.

\vspace{0.5cm}
For $j=2$, we have the following long exact sequence in homology:
\begin{center}

$0\longrightarrow \bold{S}^{(1^3)} \longrightarrow H_{1,2}(K_2,w) \longrightarrow 0 \cdots \longrightarrow 0$, 
\end{center}
from which we can conclude that $H_{1,2}(K_2,w) = \bold{S}^{(1^3)}$ and $H_{0,2}(K_2,w)= 0$.

\end{esempio}

Now, given a graph $(G,w)$, let $span_0(G,w)$ denote the homological span of the degree 0 weighted chromatic symmetric homology of $(G,w)$, i.e. of $H_{i,0}(G,w)$. We have $span_0(G,w) = k + 1$ where $k$ is maximal among indices such that $H_{k,0}(G,w)\neq 0$, {since we are assuming that $G$ has no loops, so $H_{0,0}(G,w)$ is always nonzero.}

In \cite{2}, the authors left open the following 

\begin{conjecture*}[C.6]
Let $G$ be a graph with $n$ vertices and $m$ edges, and let $b$ denote the number of blocks of $G$. Then $n-b\leq span_0(G)\leq n-1$.
\end{conjecture*}

{
We denote by $k_{max}^j(G,w)$ the largest index $k$ such that $H_{k,j}(G,w)\neq 0$ and by $k_{min}^j(G,w)$ the smallest one. 
As observed earlier, $k_{min}^0(G,w)$ is always 0. }

Using the deletion-contraction long exact sequence for weighted chromatic symmetric homology \ref{eqn} we can prove that

\begin{teo}
Let $(G,w)$ be a graph with $n$ vertices and $m$ edges. Then $k_{max}^j(G,w)\leq n-1$ for all $j\geq 0$. Moreover, if $m\geq 1$, $k_{max}^0(G,w)\leq n-2$, so $span_0(G)\leq n-1$.
\end{teo}

\begin{proof}
We prove that, if $i \geq 0$ is an index such that $H_{i,j}(G,w)\neq 0$, then we have $i \leq n - 1$. 

We proceed by induction on the number $m \geq 0$ of edges of $G$.
{If $m=0$, we have that the homology $H_{i,j}(G,w)$ is trivial for all $i>0$, since we don't have any states with more than zero edges. Therefore, the first inequality holds.}

{Furthermore, if we require $m\geq 1$, at the base step we have to consider the case $m=1$. It follow from Remark \ref{r9} that we can assume without loss of generality that $G$ is connected,} so, if $m=1$, then $G$ is a segment with two vertices and an edge between them. {It follows from Example \ref{ex1} that $k_{max}^0(G,w)=0$, so the second part of the theorem holds. }

{We} now assume the {statement} true for any graph with $m-1$ edges. Let $v(G)$ denote the number of vertices of $G$ and $e(G)$ the number of edges of $G$. We have that $v(G\setminus e)=v(G)$ and $e(G \setminus e)=e(G)-1=m-1$. Moreover, we have that $v(G / e)=v(G)-1$ and $e(G / e)=e(G)-1=m-1$. 

Let $i > v(G)-2$. Since $v(G\setminus e)=v(G)$, we have also that $i > v(G\setminus e)-2$. By inductive hypothesis, we have $H_{i,j}(G\setminus e,w)=0$. Moreover, since $i-1 > v(G) -3 = v(G/e) -2$, by inductive hypothesis, we have $H_{i-1,j}(G / e, w/e)=0$ and $H_{i,j}(G / e, w/e)=0$.

From the deletion-contraction long exact sequence \ref{eqn}
\begin{center}
  $\dots \longrightarrow H_{i,j}(G/e,w/e) \longrightarrow H_{i,j}(G\setminus e,w) \rightarrow H_{i,j}(G,w) \rightarrow  H_{i-1,j}(G/e,w/e) \longrightarrow $,  
\end{center}
it follows that $H_{i,j}(G,w)=0$.
\end{proof}



In \cite{2}, the authors left open also the following

\begin{conjecture*}[C.5]
Given any graph $G$, chromatic symmetric homology groups $H_{i,0}(G;\mathbb{C})$ are non-trivial for all $0 \leq i \leq span_0(G) - 1$, $j\geq 0$.
\end{conjecture*}

Using the deletion-contraction long exact sequence, we can prove the following

\begin{teo}\label{t7}
Let $(G,w)$ be a graph. Then $H_{i,j}(G,w;\mathbb{C})$ is non-trivial for all $k_{min}^j(G,w) \leq i \leq k_{max}^j(G,w)$, $j\geq 0$.
\end{teo}

Since $k_{min}^0(G,w)$ is always 0, Theorem \ref{t7} shows in particular that Conjecture C.5 is true.

\begin{proof}
We proceed by induction on the number $m \geq 0$ of edges of $G$. 
{If $m=0$, we have that the homology $H_{i,j}(G,w)$ is trivial for all $i>0$, since we don't have any states with more than zero edges. Therefore, the result holds.} 

Now assume the {statement} true for any graph with $m-1$ edges.

If $k_{max}^j(G \setminus e,w) \geq k_{max}^j(G,w)$, since $G \setminus e$ has $m-1$ edges, by inductive hypothesis, we have that $H_{k_{max}^j(G,w),j}(G \setminus e,w) \neq 0$. If $H_{k_{max}^j(G,w)-1,j}(G/e,w/e)=0$, then by inductive hypothesis, it is also $H_{k_{max}^j(G,w),j}(G/e,w/e)=0$. Therefore, by the deletion-contraction long exact sequence \ref{eqn}
\begin{center}
  $\longrightarrow H_{k_{max}^j(G,w),j}(G/e,w/e) \longrightarrow H_{k_{max}^j(G,w),j}(G\setminus e,w) \rightarrow H_{k_{max}^j(G,w),j}(G,w) \rightarrow  H_{k_{max}^j(G,w)-1,j}(G/e,w/e) \longrightarrow \dots$,  
\end{center}
we have $H_{k_{max}^j(G,w),j}(G\setminus e,w) \cong H_{k_{max}^j(G,w),j}(G,w)$.

Otherwise, {$H_{k_{max}^j(G,w)-1,j}(G/e,w/e)\neq0$, so} $k_{max}^j(G/e,w/e)\geq k_{max}^j(G,w)-1$.

If instead $k_{max}^j(G \setminus e,w) < k_{max}^j(G,w)$, we have $H_{k_{max}^j(G,w),j}(G\setminus e,w)=0$ and $H_{k_{max}^j(G,w),j}(G,w) \neq 0$. Therefore, by the deletion-contraction long exact sequence \ref{eqn}
\begin{center}
  $\dots \longrightarrow H_{k_{max}^j(G,w),j}(G\setminus e,w) \rightarrow H_{k_{max}^j(G,w),j}(G,w) \rightarrow  H_{k_{max}^j(G,w)-1,j}(G/e,w/e) \longrightarrow \dots$,  
\end{center}
we have that the map from $H_{k_{max}^j(G,w),j}(G,w)$ to  $H_{k_{max}^j(G,w)-1,j}(G/e,w/e)$ is injective. Hence, $H_{k_{max}^j(G,w),j}(G,w)$ is isomorphic to the image of this map, which is a non-trivial submodule of $H_{k_{max}^j(G,w)-1,j}(G/e,w/e)$. It follows that 

$H_{k_{max}^j(G,w)-1,j}(G/e,w/e) \neq 0$ and $k_{max}^j(G/e,w/e)\geq k_{max}^j(G,w)-1$.

Now assume $k_{min}^j(G,w)\leq i\leq k_{max}^j(G,w)$ and prove that $H_{i,j}(G,w)$ is non-trivial. As observed above, we have three cases to consider: 
\begin{itemize}
    \item[(i)] $k_{max}^j(G \setminus e,w) \geq k_{max}^j(G,w)$ and $H_{k_{max}^j(G,w),j}(G\setminus e,w) \cong H_{k_{max}^j(G,w),j}(G,w)$;
    \item[(ii)] $k_{max}^j(G \setminus e,w) \geq k_{max}^j(G,w)$ and $k_{max}^j(G/e,w/e)\geq k_{max}^j(G,w)-1$;
    \item[(iii)] $k_{max}^j(G \setminus e,w) < k_{max}^j(G,w)$ {and $k_{max}^j(G/e,w/e)\geq k_{max}^j(G,w)-1$}.
\end{itemize}

\hspace{-0.5cm}In case (i), $k_{max}^j(G \setminus e,w) \geq k_{max}^j(G,w)$ and $H_{k_{max}^j(G,w),j}(G\setminus e) \cong H_{k_{max}^j(G,w),j}(G,w)$, so by inductive hypothesis we have that $H_{i,j}(G \setminus e,w)$ is non-trivial. It follows {from \ref{eqn}, and for how the maps are defined,} that also $H_{i,j}(G,w)$ is non-trivial.

In case (ii), if $k_{max}^j(G \setminus e,w) \geq k_{max}^j(G,w)$ and $k_{max}^j(G/e,w/e)\geq k_{max}^j(G,w)-1$, then $i-1 \leq k_{max}^j(G,w) -1 \leq k_{max}^j(G/e,w/e)$. Therefore, by induction, $H_{i-1,j}(G/e)$ is non-trivial. Moreover, by induction, also $H_{i,j}(G \setminus e,w)$ is non trivial. It follows {from \ref{eqn}, and for how the maps are defined,} that also $H_{i,j}(G,w)$ is non-trivial.

Finally, we consider the case (iii) with $k_{max}^j(G \setminus e,w) < k_{max}^j(G,w)$. We just have to see what happens if $k_{max}^j(G \setminus e,w) < i \leq k_{max}^j(G,w)$, since, if $i \leq k_{max}^j(G \setminus e,w) < k_{max}^j(G,w)$, {as in the previous case, both $H_{i-1,j}(G/e,w/e)$ and $H_{i,j}(G \setminus e,w)$ are non-trivial, and so it is $H_{i,j}(G,w)$ $\neq$ 0}. If  $k_{max}^j(G \setminus e,w) < i \leq k_{max}^j(G,w)$, we have that $H_{i,j}(G \setminus e,w) = 0$. From the deletion-contraction long exact sequence \ref{eqn}
\begin{center}
  $\dots \longrightarrow H_{i,j}(G\setminus e,w) \rightarrow H_{i,j}(G,w) \rightarrow  H_{i-1,j}(G/e,w/e) \longrightarrow \dots$,  
\end{center}
it follows that the map from $H_{i,j}(G,w)$ to $H_{i-1,j}(G/e,w)$ is injective. Moreover, since $i-1 \leq k_{max}^j(G,w) -1 \leq k_{max}^j(G/e,w/e)$, as proved above, by induction, $H_{i-1,j}(G/e,w/e)$ is non-trivial. Hence, for how the maps are defined, $H_{i,j}(G,w)$ is non-trivial.
\end{proof}

\subsection{Future directions}
{Chandler, Sazdanovic, Stella and Yip in \cite{2} investigated the properties of chromatic symmetric homology with integer coefficients. They conjectured  that a graph $G$ is non-planar if and only if its chromatic symmetric homology in bidegree (1,0) contains $\mathbb{Z}_2$-torsion. In \cite{MR4540909}, the authors showed that the chromatic symmetric homology of a finite non-planar graph contains $\mathbb{Z}_2$-torsion in bidegree (1,0). We hope that these new tools will help to understand if this conjecture is true also in the other direction. }

{Moreover, we think that the deletion-contraction long exact sequence could simplify the computation of the homology, even in the unweighted case, and allow to study it better.}

\section*{Acknowledgments}
I thank Salvatore Stella and Luca Moci for suggesting me to work on this topic and for many helpful conversations about it; without them this paper would not have been possible. I thank Alex Chandler for reading the article and for his valuable advise. Finally, I am grateful to the reviewers for their precise and useful comments.
\printbibliography[heading=bibintoc]

@article {3,
    AUTHOR = {Stanley, Richard P.},
     TITLE = {A symmetric function generalization of the chromatic
              polynomial of a graph},
   JOURNAL = {Adv. Math.},
  FJOURNAL = {Advances in Mathematics},
    VOLUME = {111},
      YEAR = {1995},
    NUMBER = {1},
     PAGES = {166--194},
      ISSN = {0001-8708},
   MRCLASS = {05E05 (05C15)},
  MRNUMBER = {1317387},
MRREVIEWER = {SeungKyung Park},
       DOI = {10.1006/aima.1995.1020},
       URL = {https://doi.org/10.1006/aima.1995.1020},
}

@article {1,
    AUTHOR = {Sazdanovic, Radmila and Yip, Martha},
     TITLE = {A categorification of the chromatic symmetric function},
   JOURNAL = {J. Combin. Theory Ser. A},
  FJOURNAL = {Journal of Combinatorial Theory. Series A},
    VOLUME = {154},
      YEAR = {2018},
     PAGES = {218--246},
      ISSN = {0097-3165},
   MRCLASS = {05E05 (05C15 05C31)},
  MRNUMBER = {3718066},
MRREVIEWER = {Elizabeth M. Niese},
       DOI = {10.1016/j.jcta.2017.08.014},
       URL = {https://doi.org/10.1016/j.jcta.2017.08.014},
}

@book {4,
    AUTHOR = {Fulton, William},
     TITLE = {Young tableaux},
    SERIES = {London Mathematical Society Student Texts},
    VOLUME = {35},
      NOTE = {With applications to representation theory and geometry},
 PUBLISHER = {Cambridge University Press, Cambridge},
      YEAR = {1997},
     PAGES = {x+260},
      ISBN = {0-521-56144-2; 0-521-56724-6},
   MRCLASS = {05E10 (05E05 05E15 14M15 20G05)},
  MRNUMBER = {1464693},
MRREVIEWER = {Tadeusz J\'{o}zefiak},
}

@article {5,
    AUTHOR = {Bar-Natan, Dror},
     TITLE = {On {K}hovanov's categorification of the {J}ones polynomial},
   JOURNAL = {Algebr. Geom. Topol.},
  FJOURNAL = {Algebraic \& Geometric Topology},
    VOLUME = {2},
      YEAR = {2002},
     PAGES = {337--370},
      ISSN = {1472-2747},
   MRCLASS = {57M27},
  MRNUMBER = {1917056},
MRREVIEWER = {Jacob Andrew Rasmussen},
       DOI = {10.2140/agt.2002.2.337},
       URL = {https://doi.org/10.2140/agt.2002.2.337},
}

@article {MR4093019,
    AUTHOR = {Crew, Logan and Spirkl, Sophie},
     TITLE = {A deletion-contraction relation for the chromatic symmetric
              function},
   JOURNAL = {European J. Combin.},
  FJOURNAL = {European Journal of Combinatorics},
    VOLUME = {89},
      YEAR = {2020},
     PAGES = {103143, 20},
      ISSN = {0195-6698},
   MRCLASS = {05E05 (05C15 05C31)},
  MRNUMBER = {4093019},
MRREVIEWER = {Hao Qi},
       DOI = {10.1016/j.ejc.2020.103143},
       URL = {https://doi.org/10.1016/j.ejc.2020.103143},
}

@article {MR4540909,
    AUTHOR = {Ciliberti, Azzurra and Moci, Luca},
     TITLE = {On {C}hromatic {S}ymmetric {H}omology and {P}lanarity of
              {G}raphs},
   JOURNAL = {Electron. J. Combin.},
  FJOURNAL = {Electronic Journal of Combinatorics},
    VOLUME = {30},
      YEAR = {2023},
    NUMBER = {1},
     PAGES = {Paper No. 1.15--},
   MRCLASS = {05C31 (05C10 05E05 20C30 55U15)},
  MRNUMBER = {4540909},
       DOI = {10.37236/11397},
       URL = {https://doi.org/10.37236/11397},
}

@misc{2,
  doi = {10.48550/ARXIV.1911.13297},
  
  url = {https://arxiv.org/abs/1911.13297},
  
  author = {Chandler, Alex and Sazdanovic, Radmila and Stella, Salvatore and Yip, Martha},
  
  title = {On the Strength of Chromatic Symmetric Homology for graphs},
  
  publisher = {arXiv},
  
  year = {2019},
  
  copyright = {arXiv.org perpetual, non-exclusive license}
}

\end{document}